\documentclass[11pt,letterpaper]{amsart}

\usepackage[utf8]{inputenc}

\usepackage{amsmath}
\usepackage{amssymb}
\usepackage{amsthm}
\usepackage{xparse}

\usepackage[pdftex]{hyperref}

\newtheorem{counter}{subcounter}[section]
\newtheorem{cor}[counter]{Corollary}
\newtheorem{lemma}[counter]{Lemma}
\newtheorem{prop}[counter]{Proposition}

\newtheorem{theorem}[counter]{Theorem}

\newtheorem{athm}{Theorem}

\theoremstyle{definition}
\newtheorem{defn}[counter]{Definition}

\newcommand{\eqnsp}[1]{\leavevmode\xleaders\hbox{\phantom{-}}\hfill\kern0pt\[#1\]}

\newcommand{\eps}{\varepsilon}

\newcommand{\hsp}{\mskip0.5\thinmuskip}
\newcommand{\msp}{\mskip0.25\thinmuskip}

\newcommand{\acdot}{\mskip0.75\thinmuskip\cdot\mskip0.75\thinmuskip}

\newcommand{\vnt}{\mathbin{\overline{\otimes}}}

\newcommand{\A}{\mathcal{A}}
\newcommand{\B}{\mathcal{B}}

\newcommand{\D}{\mathcal{D}}

\renewcommand{\H}{\mathcal{H}}
\newcommand{\I}{\mathcal{I}}
\newcommand{\J}{\mathcal{J}}
\newcommand{\K}{\mathcal{K}}
\newcommand{\M}{\mathcal{M}}
\newcommand{\N}{\mathcal{N}}
\renewcommand{\P}{\mathcal{P}}
\newcommand{\Q}{\mathcal{Q}}
\newcommand{\Rc}{\mathcal{R}}

\newcommand{\W}{\mathcal{W}}

\newcommand{\Z}{\mathcal{Z}}

\newcommand{\C}{\mathbb{C}}

\newcommand{\G}{\mathbb{G}}
\newcommand{\Nb}{\mathbb{N}}

\newcommand{\R}{\mathbb{R}}

\DeclareMathOperator{\Aut}{Aut}

\DeclareMathOperator{\Span}{Span}

\NewDocumentCommand\abs{g}{
    \lvert\IfValueTF{#1}{#1}{\,\cdot\,}\rvert
}

\NewDocumentCommand\inner{gg}{
    \langle\IfValueTF{#2}{#1,#2}{\cdot,\IfValueTF{#1}{#1}{\cdot}}\rangle
}

\NewDocumentCommand\norm{g}{
    \lVert\IfValueTF{#1}{#1}{\,\cdot\,}\rVert
}

\NewDocumentCommand\cbnorm{g}{
    \norm{#1}_\mathrm{cb}
}

\NewDocumentCommand\opnorm{g}{
    \norm{#1}_\mathrm{op}
}

\RenewDocumentCommand\Im{s}{
    \operatorname{Im}\IfBooleanTF{#1}{\,}{\hsp}
}

\RenewDocumentCommand\Re{s}{
    \operatorname{Re}\IfBooleanTF{#1}{\,}{\hsp}
}

\calclayout

\hypersetup{
    colorlinks = true,
    linkcolor = blue,
    citecolor = red,
    urlcolor = cyan
}

\begin{document}

\title{Representation stability for compact and discrete quantum groups}

\author{Michael Brannan}
\address{University of Waterloo, 200 University Avenue West, Waterloo, Ontario, N2L 3G1, Canada}
\email{michael.brannan@uwaterloo.ca}

\author{Junichiro Matsuda}
\address{Research Institute of Mathematics, Seoul National University, 1 Gwanak-ro, Gwanak-gu, Seoul 08826, Republic of Korea \emph{and} BKC Research Organization of Social Sciences, Ritsumeikan University, 1-1-1 Noji-higashi, Kusatsu, Shiga 525-8577, Japan}
\email{junichiro.matsuda@snu.ac.kr}

\author{Erik S\'eguin}
\address{University of Waterloo, 200 University Avenue West, Waterloo, Ontario, N2L 3G1, Canada}
\email{erik.seguin@uwaterloo.ca}

\dedicatory{Dedicated to the memory of Zhong-Jin Ruan}

\thanks{The authors were partially supported by an NSERC Discovery Grant. The second author was supported by JSPS KAKENHI Grant Number JP23KJ1270.}

\begin{abstract}
We study approximate representations of locally compact quantum groups and prove stability results in the sense of Ulam in this context. Our main result is that compact and amenable discrete quantum groups are representation stable. We also show that an analogous stability result holds for unitary compressions of general amenable locally compact quantum groups without the assumption of compactness or discreteness. 
\end{abstract}

\maketitle

\setcounter{tocdepth}{1}
\tableofcontents

\section{Introduction}

In \cite{ulam60}, Ulam put forth the following question: under what circumstances do small changes in the initial hypotheses of a mathematical statement result in the conclusion still being approximately satisfied? This phenomenon has since been intensely studied in a variety of contexts, with that of group representations prominent among them. An \emph{approximate representation} of a group $G$ on a Hilbert space $\H$ is a map ${\varphi:G\to B(\H)}$ that satisfies the usual homomorphism condition only approximately in norm rather than exactly. In this context, Ulam's question may be phrased as follows: when is an approximate representation well-approximated by a genuine representation? 

Ulam stability of unitary representations is well-studied in the setting of classical locally compact groups. Stability in this sense was initially demonstrated for compact Lie groups in \cite{gkr74}; this was later extended to the case of general compact topological groups in \cite{hk77}. It was subsequently shown in \cite{kazhdan82} that amenable discrete groups are representation stable, and the general amenable locally compact case was completed a few years later as a special case of the analogous results in \cite{johnson88} for amenable Banach algebras. 

In this paper, we consider the problem of Ulam stability for unitary representations of locally compact quantum groups. We recall that a locally compact quantum group is an object $\G$ given by (among other things) a von Neumann algebra ${L^\infty(\G)}$ along with a co-associative unital normal *-homomorphism ${\Delta:L^\infty(\G)\to L^\infty(\G)\vnt L^\infty(\G)}$. The pre-adjoint of $\Delta$ endows the predual ${L^1(\G)}$ of ${L^\infty(\G)}$ with the structure of a completely contractive Banach algebra. This algebra also admits a densely defined involution $\sharp$ whose domain is denoted by ${\smash{L_\sharp^1}(\G)}$. The locally compact quantum groups for which ${L^\infty(\G)}$ is commutative are in canonical one-to-one correspondence with classical locally compact groups. Unitary representations of a locally compact quantum group $\G$ valued in a von Neumann algebra $\M$ are unitary elements ${X\in U(\M\vnt L^\infty(\G))}$ satisfying 
\eqnsp{(\iota\otimes\Delta)(X)=X_{12}\hsp X_{13}}
As in the setting of classical locally compact groups, unitary representations ${X\in U(\M\vnt L^\infty(\G))}$ are in one-to-one correspondence with completely bounded homomorphisms ${\varphi:L^1(\G)\to\M}$ which restrict to non-degenerate *-homomorphisms on ${\smash{L_\sharp^1}(\G)}$. Given ${\eps\geq0}$, one can relax the notion of a representation to that of an \emph{$\eps$-representation}: this is an element ${X\in\M\vnt L^\infty(\G)}$ such that 
\eqnsp{\norm{(\iota\otimes\Delta)(X)-X_{12}\hsp X_{13}}\leq\eps}

The following is our main result. 

\begin{athm}\label{cqg-adqg-stability}
Let ${\eps\in[0,\frac{1}{2})}$, $\G$ be a compact quantum group or amenable discrete quantum group, and $\M$ be a von Neumann algebra. If ${X\in U(\M\vnt L^\infty(\G))}$ is an $\eps$-representation, then there exists a unitary representation ${U\in U(\M\vnt L^\infty(\G))}$ such that ${\norm{X-U}\leq\eps+4\eps^2}$. 
\end{athm}

In proving the above, it is natural to consider modifying the techniques of \cite{johnson88} to the setting of completely contractive Banach algebras to obtain stability on the ${L^1}$-level. However, a particular issue presents itself with this approach: namely, that it requires that we restrict our attention to those locally compact quantum groups $\G$ for which ${L^1(\G)}$ is operator amenable. As there are many interesting examples of amenable locally compact quantum groups for which this is not the case (see e.g., \cite{clr15}), this makes for a rather displeasing constraint. We therefore opt to instead adapt the method of \cite{dcot19}, which consists of first upgrading the initial approximate representation to one which is positive definite with a certain ``improving operator'' (cf. \cite[Theorem 1.3]{shtern99}, \cite[Theorem 3.1]{bot13}) and then using a Stinespring-type dilation to obtain a unitary representation which (upon taking an appropriate cutdown) is close in norm to the original approximate representation. The intuition behind the above ``improving operator'' approach in the classical setting is that if $G$ is an amenable locally compact group with right invariant mean ${\eta\in L^\infty(G)^*}$ and ${X\in\M\vnt L^\infty(G)}$, then one can try to make sense of the $\M$-valued function on $G$ defined by the formal ``integral'' specified by 
\eqnsp{s\to\int_G X(t)^*\hsp X(ts)\,d\eta(t),\qquad s\in G}
When $G$ is compact or discrete, one can make precise sense of this $\M$-valued function and it can further be shown that it is positive definite and therefore admits a Stinespring-type dilation to a continuous unitary representation. Furthermore, if $X$ is a unitary $\eps$-representation, this ``averaged'' version of $X$ is a ${2\eps}$-representation. We show that these ideas can be made sense of in the broader setting of compact and discrete quantum groups. (We mention here the recent works \cite{at26} and \cite{dbt26}, which treat this idea in the setting of classical compact groups.)

As a consequence of this approach, we obtain a partial stability result for the general amenable locally compact case in the course of proving the above result. In the following theorem statement, a \emph{compression} refers to a corner of a unitary representation (see Definition \ref{stinespring-type-defn}). 

\begin{athm}\label{amenable-pos-def-stability-thm}
Let ${\eps\in[0,\frac{1}{4})}$, $\G$ be an amenable locally compact quantum group, and $\M$ be a von Neumann algebra. If ${X\in\M\vnt L^\infty(\G)}$ is a compression which is $\eps$-unitary in the sense that
\eqnsp{\norm{1_{\M\,\vnt\,L^\infty(\G)}-X^*X}\leq\eps,\qquad\qquad\norm{1_{\M\,\vnt\,L^\infty(\G)}-XX^*}\leq\eps}
then there exists a unitary representation ${U\in U(\M\vnt L^\infty(\G))}$ such that ${\norm{X-U}\leq4\eps}$. 
\end{athm}

The paper is organized as follows. In Section \ref{prelim-section}, we introduce our notation and recall the various necessary preliminaries. In Section \ref{cqg-section}, we prove a Stinespring-type dilation theorem for a class of improving operators on compact quantum groups. In Section \ref{adqg-section}, we prove a pointwise factorization result for the corresponding completely bounded linear maps on ${\smash{L_\sharp^1}(\G)}$ and use this to obtain an analogous dilation result in the setting of amenable discrete quantum groups. Finally, in Section \ref{stability-section} we prove our desired stability results. 
\medskip

\noindent\textbf{Acknowledgments.} The authors would like to thank Nico Spronk and Jennifer Zhu for valuable discussions on the topic of this paper. 

\section{Preliminaries}\label{prelim-section}

\subsection{Notation}

Throughout this paper, we let $\iota$ denote the identity map from a space to itself; abusing notation, we omit any reference to the underlying space, trusting that context will suffice to preclude any ambiguity in this regard. We let $\otimes$ denote the algebraic tensor product of vector spaces and $\vnt$ denote the spatial tensor product of von Neumann algebras. We employ the usual leg numbering notation on spatial tensor products (see e.g., \cite{bs93}). 

If $\H$ is a Hilbert space, then we let ${B(\H)}$ denote the space of bounded operators on $\H$ equipped with the operator norm ${\opnorm}$. If ${\xi,\eta\in\H}$, then we let ${\omega_{\xi,\eta}\in B(\H)^*}$ denote the bounded linear functional defined by 
\eqnsp{\omega_{\xi,\eta}=\inner{\acdot\hsp\xi}{\eta}}

If $\A$ is a *-algebra, then we let ${\Z(\A)}$ denote the center of $\A$ and let 
\eqnsp{\A_+=\{x^*x:x\in\A\}}
Moreover, if $\A$ is unital, then we let ${U(\A)}$ denote the unitary group of $\A$. 

If $\M$ is a von Neumann algebra, then we let ${\M_*}$ denote the predual of $\M$ and let ${\M_{*,+}}$ denote the cone of positive linear functionals in ${\M_*}$. If $\M$ is concretely represented on a Hilbert space $\H$, then we let ${\M'}$ denote the commutant of $\M$ in ${B(\H)}$. 

Furthermore, if ${\varphi:\M_+\to[0,\infty]}$ is a weight, then we let ${(\H_\varphi,\pi_\varphi,\Lambda_\varphi)}$ denote the corresponding GNS construction and let 
\eqnsp{\J_\varphi=\{x\in\M:\varphi(x^*x)<\infty\},\qquad\qquad\M_\varphi=\Span(\{y^*x:x,y\in\J_\varphi\})}
There is a unique linear functional ${\psi:\M_\varphi\to\C}$ extending ${\varphi\rvert_{(\M_\varphi)_+}}$; abusing notation, we write $\varphi$ for this linear functional as well. 

If $X$, $Y$, $Z$, and $W$ are vector spaces and ${\varphi:X\to Z}$ and ${\psi:Y\to W}$ are linear maps, then we let ${\varphi\otimes\psi:X\otimes Y\to Z\otimes W}$ denote the \emph{algebraic tensor product} of $\varphi$ and $\psi$, i.e., the unique linear map satisfying the equality 
\eqnsp{(\varphi\otimes\psi)(x\otimes y)=\varphi(x)\otimes\psi(y)}
for all ${x\in X}$ and ${y\in Y}$. 

If $\A$ and $\B$ are *-algebras and ${\varphi:\A\to\B}$ is a linear map, then we let ${\varphi^*:\A\to\B}$ denote the linear map defined by 
\eqnsp{\varphi^*(x)=\varphi(x^*)^*}

If $E$ and $F$ are operator spaces, then we let ${CB(E,F)}$ denote the space of completely bounded linear maps from $E$ to $F$. We usually write simply ${CB(E)}$ rather than ${CB(E,E)}$. 

If $\A$ is a C*-algebra and $E$ is a pre-Hilbert $\A$-module, then we let ${\mathcal{L}(E)}$ denote the space of adjointable operators on $E$. 

\subsection{Tensor amplifications and slice maps}

Throughout this subsection, we let $\M$, $\N$, $\P$, and $\Q$ denote von Neumann algebras. 

If ${\varphi:\M\to\P}$ and ${\psi:\N\to\Q}$ are completely positive normal linear maps, then (in an abusive overloading of notation) we let ${\varphi\otimes\psi:\M\vnt\N\to\P\vnt\Q}$ denote the unique (completely positive) normal linear map agreeing with the algebraic tensor product of $\varphi$ and $\psi$ on ${\M\otimes\N}$. 

We fix \emph{amplifications} ${CB(\M)\to CB(\M\vnt\N)}$ and ${CB(\M)\to CB(\N\vnt\M)}$ in the sense of \cite{neufang04}; this assignment will hold throughout the paper. If ${\varphi\in CB(\M)}$, then (overloading notation) we let ${\varphi\otimes\iota\in CB(\M\vnt\N)}$ and ${\iota\otimes\varphi\in CB(\N\vnt\M)}$ denote the images of $\varphi$ under the first and second amplification map, respectively. 

If ${\varphi\in CB(\M)}$ and ${\psi\in CB(\N)}$ are completely bounded linear maps at least one of which is normal, then 
\eqnsp{(\varphi\otimes\iota)\circ(\iota\otimes\psi)=(\iota\otimes\psi)\circ(\varphi\otimes\iota)}
and in this case (again overloading notation) we let ${\varphi\otimes\psi}$ denote the above expression. 

If ${\varphi:\M_+\to[0,\infty]}$ and ${\psi:\N_+\to[0,\infty]}$ are normal semi-finite weights, then we write ${\varphi\otimes\psi}$ for the associated \emph{tensor product weight} and we let ${\Lambda_\varphi\otimes\Lambda_\psi}$ denote the canonical GNS map for ${\varphi\otimes\psi}$. 
We refer the reader to \cite{vaes01} for details regarding these constructions. 

For every ${X\in\M\vnt\N}$, we let ${\Phi_X:\N_*\to\M}$ denote the associated evaluation on slice maps, i.e., the completely bounded linear map defined by 
\eqnsp{\Phi_X(\omega)=(\iota\otimes\omega)(X)}
The map ${X\to\Phi_X}$ gives a complete isometric isomorphism between ${\M\vnt\N}$ and ${CB(\N_*,\M)}$. 

\subsection{Locally compact quantum groups}

In this subsection, we summarize the basic background pertaining to the theory of locally compact quantum groups that we shall use throughout the paper; our exposition is minimal and the reader is referred to \cite{kv00} and \cite{kv03} for more detail. We recall that a \emph{locally compact quantum group} (or \emph{quantum group}) is a quadruple ${\G=(\M,\Delta,\varphi,\psi)}$ consisting of a von Neumann algebra $\M$, a unital normal *-homomorphism ${\Delta:\M\to\M\vnt\M}$ such that 
\eqnsp{(\Delta\otimes\iota)\circ\Delta=(\iota\otimes\Delta)\circ\Delta}
and two faithful normal semi-finite weights ${\varphi:\M_+\to[0,\infty]}$ and ${\psi:\M_+\to[0,\infty]}$ (which we refer to as the \emph{left Haar weight} and \emph{right Haar weight}, respectively) such that 
\eqnsp{\varphi((\omega\otimes\iota)(\Delta(x)))=\varphi(x)\,\omega(1_\M),\qquad\qquad\psi((\iota\otimes\omega)(\Delta(y)))=\psi(y)\,\omega(1_\M)}
for all ${x\in(\M_\varphi)_+}$, ${y\in(\M_\psi)_+}$, and ${\omega\in\M_{*,+}}$. 

We let ${L^\infty(\G)}$ denote the von Neumann algebra $\M$ and let ${L^1(\G)}$ denote its predual ${\M_*}$, and we endow the latter with the structure of a Banach algebra by equipping it with the convolution product defined by 
\eqnsp{\omega\star\eta=(\omega\otimes\eta)\circ\Delta}

There exists a unique (ultrastrong-*) densely defined closed linear map ${S:\D(S)\to L^\infty(\G)}$ (the \emph{antipode}) characterized by the following conditions: 
\begin{enumerate}
\item ${(\iota\otimes\varphi)(\Delta(x^*)\msp(1_{L^\infty(\G)}\otimes y))\in\D(S)}$ for all ${x,y\in\J_\varphi}$, and 
\item ${S((\iota\otimes\varphi)(\Delta(x^*)\msp(1_{L^\infty(\G)}\otimes y)))=(\iota\otimes\varphi)((1_{L^\infty(\G)}\otimes x^*)\,\Delta(y))}$ for all ${x,y\in\J_\varphi}$. 
\end{enumerate}
The antipode admits a ``polar decomposition'' in the following sense: there exists a unique strongly continuous one-parameter automorphism group ${\{\tau_t\}_{t\,\in\,\R}\subseteq\Aut(L^\infty(\G))}$ (the \emph{scaling group}) and a unique involutive *-anti-automorphism ${R:L^\infty(\G)\to L^\infty(\G)}$ (the \emph{unitary antipode}) satisfying the following properties: 
\begin{enumerate}
\item ${\tau_t\circ R=R\circ\tau_t}$ for all ${t\in\R}$, and 
\item ${S=R\circ\tau_{-i/2}}$. 
\end{enumerate}
There is a unique real number ${\nu>0}$ (the \emph{scaling constant}) such that ${\varphi\circ\tau_t=\nu^{-t}\hsp\varphi}$ for all ${t\in\R}$. 

Furthermore, there is a unique unitary operator ${W\in B(L^2(\G))\vnt L^\infty(\G)}$ (the \emph{right fundamental unitary}) such that 
\eqnsp{((\iota\otimes\pi_\psi)(W))(\Lambda_\psi(x)\otimes\Lambda_\psi(y))=(\Lambda_\psi\otimes\Lambda_\psi)(\Delta(x)\hsp(1_{L^\infty(\G)}\otimes y))}
for all ${x,y\in\J_\psi}$, where ${L^2(\G)}$ denotes the GNS Hilbert space corresponding to $\psi$. (We adopt here the conventions of \cite{bs93}, giving preference to the right Haar weight and right fundamental unitary over their left counterparts.)

We let ${\smash{L_\sharp^1}(\G)}$ denote the set of linear functionals ${\omega\in L^1(\G)}$ for which there exists ${\omega^\sharp\in L^1(\G)}$ such that 
\eqnsp{\omega^\sharp\rvert_{\D(S)}=\omega^*\circ S}
This is an ${\opnorm}$-norm dense subalgebra of ${L^1(\G)}$, and equipping it with the involution ${\omega\to\omega^\sharp}$ and the norm ${\norm_\sharp}$ defined by 
\eqnsp{\norm{\omega}_\sharp=\max\{\opnorm{\omega},\opnorm{\omega^\sharp}\}}
endows it with the structure of a Banach *-algebra. 

Finally, we fix some additional notation. If $\N$ is a von Neumann algebra and ${X\in\N\vnt L^\infty(\G)}$, then we let ${\Psi_X}$ denote the restriction of ${\Phi_X}$ to the subspace ${\smash{L_\sharp^1}(\G)}$.

We let ${\mu:L^\infty(\G)\otimes L^\infty(\G)\to L^\infty(\G)}$ denote the linear multiplication map, and for every ${n\in\Nb}$, we let ${\mu^{(n)}:L^\infty(\G)^{\otimes\msp(n+1)}\to L^\infty(\G)}$ denote its multilinear extension: 
\eqnsp{\mu^{(1)}=\mu,\qquad\qquad\mu^{(n+1)}=\mu^{(n)}\circ(\iota^{\otimes\msp n}\otimes\mu)}
Similarly, for every ${n\in\Nb}$, we let ${\Delta^{(n)}:L^\infty(\G)\to L^\infty(\G)^{\vnt\msp(n+1)}}$ denote the iterated co-product: 
\eqnsp{\Delta^{(1)}=\Delta,\qquad\qquad\Delta^{(n+1)}=(\iota^{\otimes\msp n}\otimes\Delta)\circ\Delta^{(n)}}

\subsection{Compact and discrete quantum groups}

A locally compact quantum group $\G$ is \emph{compact} if the left Haar weight (equivalently, the right Haar weight) is finite. In this case, the Haar weights are equal up to rescaling; henceforth, we implicitly assume that $\varphi$ and $\psi$ are both unital (therefore equal) whenever $\G$ is compact. 

A locally compact quantum group $\G$ is \emph{discrete} if ${L^1(\G)}$ is unital. In this case, there exists a collection ${\{\M_\lambda\}_{\lambda\,\in\,\I}}$ of finite dimensional matrix algebras such that ${L^\infty(\G)\cong\prod_{\lambda\,\in\,\I}\M_\lambda}$, and we let ${C_c(\G)}$ denote the ultraweakly dense subspace of ${L^\infty(\G)}$ corresponding to the algebraic direct sum ${\bigoplus_{\lambda\,\in\,\I}\M_\lambda}$ under the above isomorphism. This is a two-sided ideal in ${L^\infty(\G)}$ and is naturally equipped with the structure of a multiplier Hopf *-algebra. Let ${\epsilon\in L^1(\G)}$ denote the multiplicative identity; then ${\epsilon\rvert_{C_c(\G)}}$ agrees with the co-unit arising from the multiplier Hopf *-algebra structure. Furthermore, ${C_c(\G)\subseteq\D(S)}$ and ${S\rvert_{C_c(\G)}}$ agrees with the antipode arising from the multiplier Hopf *-algebra structure. Consequently, (abusing notation) we write simply $\epsilon$ and $S$ for the co-unit and antipode of ${C_c(\G)}$, respectively. We refer the reader to \cite{vd94} and \cite{vd98} for more details on multiplier Hopf *-algebras. 

\subsection{Amenable quantum groups}

A locally compact quantum group $\G$ is \emph{amenable} if it admits a state ${\eta\in L^\infty(\G)^*}$ such that 
\eqnsp{\eta((\iota\otimes\omega)(\Delta(x)))=\eta(x)\,\omega(1_{L^\infty(\G)})}
for all ${x\in L^\infty(\G)}$ and ${\omega\in L^1(\G)}$. We refer to any such state as a \emph{right invariant mean}. It is clear that the normalized right Haar weight on any compact quantum group is a right invariant mean. We remark that it follows from the above equality that 
\eqnsp{\omega((\eta\otimes\iota)(\Delta(x)))=\omega(\eta(x)\,1_{L^\infty(\G)})}
for all ${x\in L^\infty(\G)}$ and ${\omega\in L^1(\G)}$, wherefore 
\eqnsp{(\eta\otimes\iota)\circ\Delta=\eta(\acdot)\,1_{L^\infty(\G)}}

\subsection{Representations and approximate representations}

Throughout this subsection, we let $\G$ denote a locally compact quantum group and $\M$ denote a von Neumann algebra. We recall that a \emph{representation} of $\G$ in $\M$ is an invertible element ${X\in\M\vnt L^\infty(\G)}$ such that 
\eqnsp{(\iota\otimes\Delta)(X)=X_{12}\hsp X_{13}}

The following result is well-known in the case where ${\M=B(\H)}$ for some Hilbert space $\H$ (see e.g., \cite{bds13}, \cite{bk18}, \cite{by19}); it is likely known to experts in the general setting as well, but we were unable to find a reference in the literature, so we provide a proof for completeness. 

\begin{prop}
The (restriction of the) map ${X\to\Psi_X}$ gives a one-to-one correspondence between unitary representations of $\G$ in $\M$ and non-degenerate *-homomorphisms from ${\smash{L_\sharp^1}(\G)}$ to $\M$. 
\end{prop}

\begin{proof}
It may be assumed without loss of generality that $\M$ and ${L^\infty(\G)}$ are concretely represented on Hilbert spaces $\H$ and $\K$, respectively. Let ${X\in\M\vnt L^\infty(\G)}$ be a unitary representation; then a straightforward adaptation of the proof of \cite[Proposition 5.2]{kustermans01} shows that ${\Psi_X}$ is a non-degenerate *-homomorphism. Now let ${\pi:\smash{L_\sharp^1}(\G)\to\M}$ be a non-degenerate *-homomorphism; then it follows from \cite[Corollary 4.3]{kustermans01} that there exists a unique unitary representation ${X\in B(\H)\vnt L^\infty(\G)}$ such that ${\pi=\Psi_X}$. As ${\smash{L_\sharp^1}(\G)}$ is dense in ${L^1(\G)}$ and ${\Phi_X(\smash{L_\sharp^1}(\G))\subseteq\M}$, it follows that ${\Phi_X(L^1(\G))\subseteq\M}$. This implies that 
\eqnsp{(\iota\otimes\omega)((S\otimes I_\K)\hsp X)=S\,\Phi_X(\omega\rvert_{L^\infty(\G)})=\Phi_X(\omega\rvert_{L^\infty(\G)})\,S=(\iota\otimes\omega)(X\hsp(S\otimes I_\K))}
for all ${\omega\in B(\K)_*}$ and ${S\in\M'}$, and so ${(S\otimes I_\K)\hsp X=X\hsp(S\otimes I_\K)}$ for all ${S\in\M'}$. It is further clear that ${(I_\H\otimes T)\hsp X=X\hsp(I_\H\otimes T)}$ for all ${T\in L^\infty(\G)'}$, whence ${X\in(\M'\vnt L^\infty(\G)')'}$; the commutation theorem for spatial tensor products then implies that ${X\in\M\vnt L^\infty(\G)}$. 
\end{proof}

Motivated by the above, we make the following definition. 

\begin{defn}
Let ${\eps\geq0}$. An element ${X\in\M\vnt L^\infty(\G)}$ is an \emph{$\eps$-representation} if 
\eqnsp{\norm{(\iota\otimes\Delta)(X)-X_{12}\hsp X_{13}}\leq\eps}
\end{defn}

It is straightforward to show (passing backward through the above one-to-one correspondence) that every bounded weakly measurable $\eps$-representation of a classical locally compact group gives rise to an $\eps$-representation in the above sense on the canonical corresponding quantum group with commutative underlying von Neumann algebra. 

\subsection{Stinespring-type elements}\label{stinespring-subsection}

A central role in this work will be played by compressions (in a particular sense) of unitary representations; we provide the precise details of our required notion below. Throughout this subsection, we let $\G$ denote a locally compact quantum group and let $\M$ denote a von Neumann algebra. 

\begin{defn}\label{stinespring-type-defn}
An element ${X\in\M\vnt L^\infty(\G)}$ is \emph{Stinespring-type} if there exists 
\begin{enumerate}
\item a von Neumann algebra $\N$,
\item a projection ${f\in\N}$,
\item a *-isomorphism ${\sigma:\M\to f\N f}$,
\item an element ${v\in\N f}$, and 
\item a unitary representation ${U\in U(\N\vnt L^\infty(\G))}$ 
\end{enumerate}
such that 
\eqnsp{(\sigma\otimes\iota)(X)=(v^*\otimes1_{L^\infty(\G)})\,U\,(v\otimes1_{L^\infty(\G)})}
\end{defn}

We refer to any quintuple ${(\N,f,\sigma,v,U)}$ witnessing the conditions of the preceding definition as a \emph{Stinespring-type dilation} of $X$. Furthermore, we say that $X$ is a \emph{compression} if it admits a Stinespring-type dilation ${(\N,f,\sigma,v,U)}$ in which $v$ is a partial isometry with ${v^*v=f}$. 

The following proposition reformulates the above notion in terms of *-representations of ${\smash{L_\sharp^1}(\G)}$. 

\begin{prop}\label{stinespring-equiv-pos-def}
Let ${X\in\M\vnt L^\infty(\G)}$. The following statements are equivalent: 
\begin{enumerate}
\item $X$ is Stinespring-type. 
\item There exists 
\begin{enumerate}
\item a von Neumann algebra $\N$,
\item a projection ${f\in\N}$,
\item a *-isomorphism ${\sigma:\M\to f\N f}$,
\item an element ${v\in\N f}$, and 
\item a non-degenerate *-homomorphism ${\pi:\smash{L_\sharp^1}(\G)\to\N}$ 
\end{enumerate}
such that 
\eqnsp{\sigma\circ\Psi_X=v^*\,\pi(\acdot)\,v}
\end{enumerate}
Furthermore, the above equivalence still holds if $X$ is assumed to be a compression in ${(1)}$ and $v$ is assumed to be a partial isometry with ${v^*v=f}$ in ${(2)}$. 
\end{prop}

\begin{proof}
Suppose that ${(1)}$ holds and let ${(\N,f,\sigma,v,U)}$ be a Stinespring-type dilation of $X$. As $U$ is a unitary representation, in turn ${\Psi_U}$ is a non-degenerate *-homomorphism; then 
\eqnsp{(\sigma\circ\Psi_X)(\omega)=\sigma((\iota\otimes\omega)(X))=(\iota\otimes\omega)((\sigma\otimes\iota)(X))=v^*\,(\iota\otimes\omega)(U)\,v=v^*\,\Psi_U(\omega)\,v}
for all ${\omega\in\smash{L_\sharp^1}(\G)}$, which yields the implication ${(1)\Rightarrow(2)}$. Now suppose instead that ${(2)}$ holds and let ${(\N,f,\sigma,v,\pi)}$ be a quintuple witnessing ${(2)}$. As $\pi$ is a non-degenerate *-homomorphism, there exists a unitary representation ${U\in U(\N\vnt L^\infty(\G))}$ such that ${\pi=\Psi_U}$. This implies that 
\eqnsp{(\iota\otimes\omega)((\sigma\otimes\iota)(X))=\sigma((\iota\otimes\omega)(X))=(\sigma\circ\Psi_X)(\omega)=v^*\,\pi(\omega)\,v=v^*\,(\iota\otimes\omega)(U)\,v}
for all ${\omega\in\smash{L_\sharp^1}(\G)}$, wherefore in turn 
\eqnsp{(\iota\otimes\omega)((\sigma\otimes\iota)(X))=v^*\,(\iota\otimes\omega)(U)\,v=(\iota\otimes\omega)((v^*\otimes1_{L^\infty(\G)})\,U\,(v\otimes1_{L^\infty(\G)}))}
for all ${\omega\in L^1(\G)}$, which yields the implication ${(2)\Rightarrow(1)}$. 
\end{proof}

It is clear from the above that ${\Psi_X}$ is completely positive if $X$ is Stinespring-type; we have the following partial converse when $\G$ is discrete. 

\begin{prop}\label{dqg-ucp-descends-compression}
Let $\G$ be discrete and ${X\in\M\vnt L^\infty(\G)}$. If ${\Psi_X}$ is unital and completely positive, then $X$ is a compression. 
\end{prop}

\begin{proof}
By \cite[Theorem 5.2]{paschke73} and the discussion following \cite[Corollary 5.3]{paschke73}, there exists a self-dual Hilbert $\M$-module $E$, a unital *-homomorphism ${\pi:\smash{L_\sharp^1}(\G)\to\mathcal{L}(E)}$, and some element ${e\in E}$ such that ${\Psi_X=\inner{\pi(\acdot)\,e}{e}}$. Let ${\N=\mathcal{L}(E)}$; then \cite[Proposition 3.10]{paschke73} implies that $\N$ is a W*-algebra. Proceeding as in the discussion after \cite[Theorem 5.2]{paschke73}, let ${f\in\N}$ be the projection defined by 
\eqnsp{f(w)=e\cdot\inner{w}{e}}
and let ${\sigma:\M\to f\N f}$ be the *-isomorphism defined by 
\eqnsp{\sigma(x)(w)=e\cdot(x\,\inner{w}{e})}
Let ${v=f}$; then ${\sigma\circ\Psi_X=v^*\,\pi(\acdot)\,v}$, and thus Proposition \ref{stinespring-equiv-pos-def} implies that $X$ is a compression. 
\end{proof}

\subsection{Miscellaneous lemmas}

We collect in this subsection a number of useful lemmas concerning general locally compact quantum groups that we shall have need of later; these results are likely known to experts, but we include them for the sake of completeness. Throughout the remainder of this subsection, we let $\G$ denote a locally compact quantum group. 

\begin{lemma}\label{antipode-central-proj}
Let ${e\in\D(S)}$ be a central projection. The following statements hold: 
\begin{enumerate}
\item ${S(e)}$ is a central projection. 
\item ${S(e)\in\D(S)}$ and ${S(S(e))=e}$. 
\end{enumerate}
\end{lemma}

\begin{proof}
It is clear that ${S(e)}$ is an idempotent. Let ${x\in\D(S)}$; then 
\eqnsp{S(e)\,S(x)=S(xe)=S(ex)=S(x)\,S(e)}
As ${S(\D(S))}$ is ultrastrong-* dense in ${L^\infty(\G)}$, it follows by the above equality that ${S(e)}$ is central. In particular, this implies that ${S(e)}$ is a normal idempotent, therefore a projection. Furthermore, as $e$ and ${S(e)}$ are self-adjoint, it follows that ${S(e)=S(e)^*\in\D(S)}$ and 
\eqnsp{S(S(e))=S(S(e)^*)=e^*=e\qedhere}
\end{proof}

\begin{lemma}\label{omega-sharp-central-supp}
Let ${\omega\in L_\sharp^1(\G)}$ and ${e,f\in\Z(L^\infty(\G))}$ be the central supports of $\omega$ and ${\omega^\sharp}$, respectively. If ${e,f\in\D(S)}$, then ${S(e)=f}$. 
\end{lemma}

\begin{proof}
Let ${\theta\in L^1(\G)}$ be the normal linear functional defined by 
\eqnsp{\theta=\omega^\sharp(S(e)\msp\acdot)}
Lemma \ref{antipode-central-proj} implies that ${S(e)\,x\in\D(S)}$ for all ${x\in\D(S)}$, and thus 
\eqnsp{\theta(x)=\omega^\sharp(S(e)\,x)=(\omega^*\circ S)(S(e)\,x)=\omega^*(S(x)\,e)=(\omega^*\circ S)(x)}
for all ${x\in\D(S)}$; it therefore follows that ${\theta\rvert_{\D(S)}=\omega^*\circ S}$, whence ${\omega^\sharp=\theta}$. As ${S(e)}$ is central, this implies that ${f\leq S(e)}$; as ${\smash{(\omega^\sharp)^\sharp}=\omega}$, the same reasoning implies that ${e\leq S(f)}$, and thus 
\eqnsp{S(e)=S(S(f)\,e)=S(e)\,f=f\qedhere}
\end{proof}

In addition to the preceding lemmas concerning projections in the domain of the antipode, we will require the following two density results. 

\begin{lemma}\label{normal-functional-dense-subspace}
If ${\I\subseteq\J_\varphi}$ is an ultraweakly dense two-sided ideal of ${L^\infty(\G)}$, then 
\eqnsp{\W=\Span(\{\varphi(y^*\cdot\,x):x,y\in\I\})}
is norm dense in ${L^1(\G)}$. 
\end{lemma}

\begin{proof}
By replacing $\I$ with its norm closure if needed, it may be assumed without loss of generality that $\I$ is norm closed; as $\I$ is a two-sided ideal in ${L^\infty(\G)}$, this implies that $\I$ is *-closed. Suppose now that $\W$ is not dense in ${L^1(\G)}$; then the Hahn-Banach separation theorem implies that there exists a nonzero element ${w\in L^\infty(\G)}$ such that ${\varphi(y^*wx)=0}$ for all ${x,y\in\I}$. As $\I$ is ultraweakly dense in ${L^\infty(\G)}$, it follows by Kaplansky's density theorem that there exists a bounded net ${(x_\alpha)_{\alpha\,\in\,\J}}$ in $\I$ such that ${(x_\alpha)\to1_{L^\infty(\G)}}$ in the ultrastrong-* topology. As ${wx_\alpha\in\I}$ for all ${\alpha\in\J}$, it follows that ${\varphi(x_\alpha^*\hsp w^*wx_\alpha)=0}$ for all ${\alpha\in\J}$; as ${(x_\alpha^*\hsp w^*wx_\alpha)\to w^*w}$ ultraweakly, the normality of $\varphi$ then implies that ${\varphi(w^*w)\leq0}$, and thus ${\varphi(w^*w)=0}$. However, as $\varphi$ is faithful, it then follows in turn that ${w=0}$, which is a contradiction; thus $\W$ is norm dense in ${L^1(\G)}$. 
\end{proof}

The proof of the following lemma is essentially contained in that of \cite[Lemma 3]{ds13}; we reproduce the details here for the convenience of the reader. 

\begin{lemma}\label{relative-operator-sharp-norm-closures}
If ${\W\subseteq\smash{L_\sharp^1}(\G)}$ is a subspace such that ${\omega\circ\tau_t\in\W}$ for all ${\omega\in\W}$ and ${t\in\R}$, then the relative $\opnorm$-norm closure of $\W$ coincides with the ${\norm_\sharp}$-norm closure of $\W$. 
\end{lemma}

\begin{proof}
It is clear that the ${\norm_\sharp}$-norm closure of $\W$ is contained in the $\opnorm$-norm closure of $\W$. For every ${n\in\Nb}$, let ${\Rc_n:L^1(\G)\to L_\sharp^1(\G)}$ be the linear map defined by 
\eqnsp{\Rc_n(\omega)=\frac{n}{\sqrt{\pi}}\int_{-\infty}^\infty e^{-n^2t^2}\hsp(\omega\circ\tau_t)\,dt}
Suppose that ${\omega\in\W}$; then ${\omega\circ\tau_t\in\W}$ and ${(\omega\circ\tau_t)^\sharp=\omega^\sharp\circ\tau_t}$ for all ${t\in\R}$. The latter implies that the map ${t\to\omega\circ\tau_t}$ is continuous with respect to the ${\norm_\sharp}$-norm topology, and thus ${\Rc_n(\omega)}$ belongs to the ${\norm_\sharp}$-norm closure of $\W$ for all ${n\in\Nb}$. Suppose now instead that ${(\omega_m)}$ is a sequence in $\W$ such that ${\opnorm{\omega_m-\omega}\to0}$ for some ${\omega\in\smash{L_\sharp^1}(\G)}$; then 
\eqnsp{\norm{\Rc_n(\omega_m)-\Rc_n(\omega)}_\sharp=\norm{\Rc_n(\omega_m-\omega)}_\sharp\leq\smash{e^{n^2/4}}\,\opnorm{\omega_m-\omega}\to0}
for all ${n\in\Nb}$. This clearly implies that ${\Rc_n(\omega)}$ belongs to the ${\norm_\sharp}$-norm closure of $\W$ for all ${n\in\Nb}$, and as ${\norm{\Rc_n(\omega)-\omega}_\sharp\to0}$, it then follows that $\omega$ belongs to the ${\norm_\sharp}$-norm closure of $\W$. 
\end{proof}

\subsection{Operator inequalities}

We will require a number of operator inequalities. All of these can be found in some form in either \cite{dcot19} or \cite{seguin22}, but we include them here for the reader's convenience. 

\begin{lemma}\label{square-ineq}
Let $\A$ be a unital C*-algebra. If ${x\in\A_+}$, then 
\eqnsp{\norm{1_\A-x}\leq\norm{1_\A-x^2}}
\end{lemma}

\begin{proof}
As ${0\leq\abs{1-t}\leq\abs{1-t^2}}$ for all ${t\in[0,\infty)}$, it follows that 
\eqnsp{\norm{1_\A-x}=\norm{\abs{1_\A-x}}\leq\norm{\abs{1_\A-x^2}}=\norm{1_\A-x^2}\qedhere}
\end{proof}

\begin{lemma}\label{expand-ineq}
Let $\A$ be a C*-algebra. If ${e\in\A}$ is a projection, ${w\in\A_+}$ with ${\norm{w}\leq1}$, and ${x\in\A e}$ with ${\norm{x}\leq1}$, then 
\eqnsp{\norm{e-x^*x}\leq\norm{e-x^*wx}}
\end{lemma}

\begin{proof}
As 
\eqnsp{0\leq x^*x=ex^*xe\leq\norm{x^*x}\,e\leq e,\qquad\qquad 0\leq x^*wx\leq\norm{w}\,x^*x\leq x^*x}
it follows that ${0\leq e-x^*x\leq e-x^*wx}$, which implies the claimed norm inequality. 
\end{proof}

\begin{lemma}\label{spectral-proj-ineq-ii}
Let $\M$ be a von Neumann algebra and ${\lambda\in[0,1)}$. If ${e\in\M}$ is a projection, ${u\in\M}$ is a partial isometry such that ${u^*u=e}$, and ${x\in\M_+}$ with ${\norm{x}\leq1}$, then 
\eqnsp{\norm{e-u^*\,\chi_{[\lambda,\msp1]}(x)\,u}\leq(1-\lambda)^{-1}\msp\norm{e-u^*xu}}
\end{lemma}

\begin{proof}
As ${0\leq1-\chi_{[\lambda,\msp1]}(t)\leq(1-\lambda)^{-1}\msp(1-t)}$ for all ${t\in[0,1]}$, it follows that  
\eqnsp{0\leq u^*\,(1_\M-\chi_{[\lambda,\msp1]}(x))\,u\leq u^*\,((1-\lambda)^{-1}\msp(1_\M-x))\,u=(1-\lambda)^{-1}\msp(e-u^*xu)}
This implies the claimed norm inequality. 
\end{proof}

\begin{lemma}\label{spectral-proj-ineq}
Let $\M$ be a von Neumann algebra and ${\lambda\in(0,1]}$. If ${w\in\M_+}$ and ${x\in\M}$, then 
\eqnsp{\norm{x^*\hsp\chi_{[\lambda,\msp1]}(w)\,x}\leq\lambda^{-1}\msp\norm{x^*wx}}
\end{lemma}

\begin{proof}
As ${0\leq\chi_{[\lambda,\msp1]}(t)\leq\lambda^{-1}\msp t}$ for all ${t\in[0,\infty)}$, in turn ${0\leq x^*\hsp\chi_{[\lambda,\msp1]}(w)\,x\leq\lambda^{-1}\msp x^*wx}$. This yields the claimed norm inequality. 
\end{proof}

\section{Compact quantum groups}\label{cqg-section}

Throughout the entirety of this section, we let $\G$ denote a compact quantum group and let $\M$ denote a von Neumann algebra. Moreover, we let ${\xi_\psi=\Lambda_\psi(1_{L^\infty(\G)})}$ and let ${E\in B(L^2(\G))}$ denote the orthogonal projection defined by 
\eqnsp{E=\inner{\acdot}{\xi_\psi}\,\xi_\psi}

\begin{lemma}\label{cqg-twist-mean-eq}
If ${X,Y\in\M\vnt L^\infty(\G)}$, then 
\eqnsp{(\iota\otimes\psi\otimes\iota)(Y_{12}^*\hsp\msp(\iota\otimes\Delta)(X))=(\iota\otimes\omega_{\xi_\psi,\hsp\xi_\psi}\otimes\iota)(((\iota\otimes\pi_\psi)(Y^*))_{12}\hsp W_{23}\hsp((\iota\otimes\pi_\psi)(X))_{12})}
\end{lemma}

\begin{proof}
As $\G$ is compact, it follows that 
\eqnsp{(\omega_{\Lambda_\psi(x),\hsp\Lambda_\psi(y)}\otimes\iota)(W)=(\psi\otimes\iota)((y^*\otimes1_{L^\infty(\G)})\,\Delta(x))}
for all ${x,y\in L^\infty(\G)}$. Let ${w\in\M}$ and ${\theta\in\M_*}$, and let ${x=(\theta\otimes\iota)((w^*\otimes1_{L^\infty(\G)})\,X)}$; then 
\begin{align*}
(\theta\otimes\psi\otimes\iota)((w\otimes y)_{12}^*\hsp\msp(\iota\otimes\Delta)(X))&=(\theta\otimes\psi\otimes\iota)((w^*\otimes y^*\otimes1_{L^\infty(\G)})\msp(\iota\otimes\Delta)(X))
\\&=(\psi\otimes\iota)((y^*\otimes1_{L^\infty(\G)})\,\Delta(x))
\\&=(\omega_{\Lambda_\psi(x),\hsp\Lambda_\psi(y)}\otimes\iota)(W)
\\&=(\omega_{\xi_\psi,\hsp\xi_\psi}\otimes\iota)((\pi_\psi(y)^*\otimes1_{L^\infty(\G)})\,W\hsp(\pi_\psi(x)\otimes 1_{L^\infty(\G)}))
\\&=(\theta\otimes\omega_{\xi_\psi,\hsp\xi_\psi}\otimes\iota)((w^*\otimes\pi_\psi(y)^*)_{12}\hsp W_{23}\hsp((\iota\otimes\pi_\psi)(X))_{12})
\\&=(\theta\otimes\omega_{\xi_\psi,\hsp\xi_\psi}\otimes\iota)(((\iota\otimes\pi_\psi)(w^*\otimes y^*))_{12}\hsp W_{23}\hsp((\iota\otimes\pi_\psi)(X))_{12})
\end{align*}
for all ${y\in L^\infty(\G)}$. Taking linear combinations and ultraweak limits then implies that 
\eqnsp{(\theta\otimes\psi\otimes\iota)(Y_{12}^*\hsp\msp(\iota\otimes\Delta)(X))=(\theta\otimes\omega_{\xi_\psi,\hsp\xi_\psi}\otimes\iota)(((\iota\otimes\pi_\psi)(Y^*))_{12}\hsp W_{23}\hsp((\iota\otimes\pi_\psi)(X))_{12})}
As the above equality holds for all ${\theta\in\M_*}$, the claimed equality follows. 
\end{proof}

We can now prove the dilation theorem alluded to in the introduction. 

\begin{theorem}\label{cqg-answ-dilation}
Let ${X,Y\in\M\vnt L^\infty(\G)}$. If ${Z\in\M\vnt L^\infty(\G)}$ is the element defined by 
\eqnsp{Z=(\iota\otimes\psi\otimes\iota)(Y_{12}^*\hsp\msp(\iota\otimes\Delta)(X))}
then 
\eqnsp{(\sigma\otimes\iota)(Z)=(v^*\otimes1_{L^\infty(\G)})\hsp W_{23}\hsp(u\otimes1_{L^\infty(\G)})}
where ${\sigma:\M\to(1_\M\otimes E)(\M\vnt B(L^2(\G)))(1_\M\otimes E)}$ is the *-isomorphism defined by 
\eqnsp{\sigma(x)=x\otimes E}
and ${u,v\in\M\vnt B(L^2(\G))}$ are the elements defined by 
\eqnsp{u=((\iota\otimes\pi_\psi)(X))(1_\M\otimes E),\qquad\qquad v=((\iota\otimes\pi_\psi)(Y))(1_\M\otimes E)}
\end{theorem}

\begin{proof}
As 
\begin{align*}
((\sigma\otimes\iota)\circ(\iota\otimes\omega_{\xi_\psi,\hsp\xi_\psi}\otimes\iota))(w\otimes x\otimes y)&=w\otimes\omega_{\xi_\psi,\hsp\xi_\psi}(x)\,E\otimes y
\\&=w\otimes ExE\otimes y
\\&=(1_\M\otimes E)_{12}\hsp\msp(w\otimes x\otimes y)\hsp(1_\M\otimes E)_{12}
\end{align*}
for all ${w\in\M}$, ${x\in B(L^2(\G))}$, and ${y\in L^\infty(\G)}$, by taking linear combinations and ultraweak limits it follows that 
\eqnsp{(\sigma\otimes\iota)\circ(\iota\otimes\omega_{\xi_\psi,\hsp\xi_\psi}\otimes\iota)=(1_\M\otimes E)_{12}\hsp\msp(\acdot)\hsp(1_\M\otimes E)_{12}}
Lemma \ref{cqg-twist-mean-eq} then clearly implies the claimed equality. 
\end{proof}

As a consequence of the above theorem, we obtain the following analogue of \cite[Proposition 2.2]{dcot19} in the setting of compact quantum groups. 

\begin{cor}\label{cqg-twist-pos-def}
Let ${X\in\M\vnt L^\infty(\G)}$. The element ${Y\in\M\vnt L^\infty(\G)}$ defined by 
\eqnsp{Y=(\iota\otimes\psi\otimes\iota)(X_{12}^*\hsp\msp(\iota\otimes\Delta)(X))}
is Stinespring-type. Furthermore, if $X$ is an isometry, then $Y$ is a compression. 
\end{cor}

\section{Amenable discrete quantum groups}\label{adqg-section}

Throughout the entirety of this section, we let $\G$ denote a discrete quantum group and let $\M$ denote a von Neumann algebra. Moreover, we let ${\smash{L_f^1}(\G)\subseteq L^1(\G)}$ be the subspace defined by 
\eqnsp{L_f^1(\G)=\Span(\{\varphi(y^*\cdot\,x):x,y\in C_c(\G)\})}

We begin this section by collecting a number of useful lemmas concerning the above subspace, which should be viewed as the set of ``finitely supported'' functionals in ${L^1(\G)}$. 

\begin{lemma}\label{central-supp-finite-dim}
If ${\omega\in\smash{L_f^1}(\G)}$, then ${\omega\in\smash{L_\sharp^1}(\G)}$ and the central support of $\omega$ belongs to ${C_c(\G)}$. 
\end{lemma}

\begin{proof}
As ${\omega\in\smash{L_f^1}(\G)}$, there exist ${x_1,\dots,x_n,y_1,\dots,y_n\in C_c(\G)}$ such that 
\eqnsp{\omega=\sum_{j\,=\,1}^n\varphi(y_j^*\cdot\msp x_j)}
Let ${e\in\Z(L^\infty(\G))}$ be the central support of $\omega$, and for ${j=1,\dots,n}$, let ${e_j\in\Z(L^\infty(\G))}$ be the central support of ${x_j}$. It is clear that ${e_j\in C_c(\G)}$ for ${j=1,\dots,n}$; as ${\bigvee_{j\,=\,1}^n e_j\in\mathrm{Alg}(\{e_1,\dots,e_n\})}$, it follows that ${\bigvee_{j\,=\,1}^n e_j\in C_c(\G)}$. As ${C_c(\G)}$ is a two-sided ideal and ${e\leq\bigvee_{j\,=\,1}^n e_j}$, this then implies that ${e\in C_c(\G)}$. Now let ${\J\subseteq C_c(\G)}$ be the compression defined by 
\eqnsp{\J=S(e)\,L^\infty(\G)}
and let ${\theta:L^\infty(\G)\to\C}$ be the linear functional defined by 
\eqnsp{\theta=(\omega^*\circ S)(S(e)\acdot\msp)}
As ${S(e)\in C_c(\G)}$, evidently $\J$ is finite dimensional. It follows that ${(\omega^*\circ S)\rvert_\J}$ is norm continuous and that the subspace topologies induced by the norm and ultraweak topologies agree on $\J$, and thus $\theta$ is normal. Moreover, it is clear that ${\theta\rvert_{\D(S)}=\omega^*\circ S}$, wherefore ${\omega\in\smash{L_\sharp^1}(\G)}$ with ${\omega^\sharp=\theta}$. 
\end{proof}

The following result is well-known to experts, but we include a short proof for completeness. 

\begin{lemma}\label{scaling-group-multiplier-hopf-inv}
If ${x\in C_c(\G)}$, then ${\tau_t(x)\in C_c(\G)}$ for all ${t\in\R}$. 
\end{lemma}

\begin{proof}
For every ${t\in\R}$, as ${\tau_t}$ is a *-isomorphism, it maps the set of minimal projections of ${L^\infty(\G)}$ onto itself; as ${C_c(\G)}$ is the linear span of the minimal projections of ${L^\infty(\G)}$, the claim follows. 
\end{proof}

The following density result will allow us to work with ${\smash{L_f^1}(\G)}$ rather than ${\smash{L_\sharp^1}(\G)}$, thereby giving us access to the playground of multiplier Hopf *-algebras via the former's connection with ${C_c(\G)}$. 

\begin{lemma}\label{finite-supp-functional-dense-subspace}
The subspace ${\smash{L_f^1}(\G)}$ is ${\norm_\sharp}$-norm dense in ${\smash{L_\sharp^1}(\G)}$. 
\end{lemma}

\begin{proof}
Let ${x_1,\dots,x_n,y_1,\dots,y_n\in C_c(\G)}$ and let ${\omega\in\smash{L_f^1}(\G)}$ be the linear functional defined by 
\eqnsp{\omega=\sum_{j\,=\,1}^n\varphi(y_j^*\cdot\msp x_j)}
For every ${t\in\R}$ and ${j=1,\dots,n}$, let ${z_{t,j}\in L^\infty(\G)}$ and ${w_{t,j}\in L^\infty(\G)}$ be the elements defined by 
\eqnsp{z_{t,j}=\tau_t(x_j),\qquad\qquad w_{t,j}=\tau_t(y_j)}
Lemma \ref{scaling-group-multiplier-hopf-inv} implies that ${z_{t,j},w_{t,j}\in C_c(\G)}$ for all ${t\in\R}$ and ${j=1,\dots,n}$; as 
\begin{align*}
\omega\circ\tau_t&=\sum_{j\,=\,1}^n\varphi(y_j^*\,\tau_t(\acdot)\,x_j)
\\&=\sum_{j\,=\,1}^n(\varphi\circ\tau_t)((w_{-t,j})^*\cdot\,z_{-t,j})
\\&=\nu^{-t}\sum_{j\,=\,1}^n\varphi((w_{-t,j})^*\cdot\,z_{-t,j})\in L_f^1(\G)
\end{align*}
for all ${t\in\R}$, it then follows by Lemma \ref{normal-functional-dense-subspace} and Lemma \ref{relative-operator-sharp-norm-closures} that ${\smash{L_f^1}(\G)}$ is dense in ${\smash{L_\sharp^1}(\G)}$. 
\end{proof}

For every ${\omega\in\smash{L_f^1}(\G)}$ and ${n\in\Nb}$, we let ${\pi_{\omega,n}:L^\infty(\G)^{\vnt\msp(n+1)}\to C_c(\G)^{\otimes\msp(n+1)}}$ denote the bounded linear map defined by 
\eqnsp{\pi_{\omega,n}(X)=(e^{\otimes\msp n}\otimes S(e))\,X}
where ${e\in\Z(L^\infty(\G))}$ is the central support of $\omega$, and we let ${\chi_{\omega,n}:L^\infty(\G)^{\vnt\msp(n+1)}\to\C}$ denote the linear functional defined by 
\eqnsp{\chi_{\omega,n}=\omega\circ\mu^{(n)}\circ(\iota^{\otimes\msp n}\otimes S)\circ\pi_{\omega,n}}

We collect here a few essential properties of the above functionals. 

\begin{lemma}\label{dqg-fin-supp-antipode-lemma-i}
If ${\omega\in\smash{L_f^1}(\G)}$ and ${n\in\Nb}$, then ${\chi_{\omega,n}\in(L^\infty(\G)^{\vnt\msp(n+1)})_*}$ and 
\eqnsp{\chi_{\omega,n}\rvert_{C_c(\G)^{\otimes(n+1)}}=\omega\circ\mu^{(n)}\circ(\iota^{\otimes\msp n}\otimes S)}
\end{lemma}

\begin{proof}
Let ${e\in\Z(L^\infty(\G))}$ be the central support of $\omega$ and ${\J\subseteq C_c(\G)^{\otimes\msp(n+1)}}$ be the compression defined by 
\eqnsp{\J=(e^{\otimes\msp n}\otimes S(e))\,L^\infty(\G)^{\vnt\msp(n+1)}}
It follows by Lemma \ref{central-supp-finite-dim} that $\J$ is finite dimensional, wherefore ${(\omega\circ\mu^{(n)}\circ(\iota^{\otimes\msp n}\otimes S))\rvert_\J}$ is norm continuous and the subspace topologies induced by the norm and ultraweak topologies agree on $\J$. It is clear that ${\pi_{\omega,n}}$ is normal; the latter of the preceding statements therefore implies that ${\pi_{\omega,n}}$ is ultraweak-to-norm topology continuous, which proves the first part of the claim. The second part of the claim is clear. 
\end{proof}

\begin{lemma}\label{dqg-fin-supp-antipode-lemma-ii}
If ${\omega\in\smash{L_f^1}(\G)}$, then ${\chi_{\omega,2}\circ(\iota\otimes\Delta)=\omega\otimes\epsilon}$. 
\end{lemma}

\begin{proof}
It follows by Lemma \ref{dqg-fin-supp-antipode-lemma-i} that ${\chi_{\omega,2}\in(L^\infty(\G)^{\vnt\msp3})_*}$; as ${\iota\otimes\Delta}$ is normal, its composition with the preceding functional is normal as well. Let ${e\in\Z(L^\infty(\G))}$ be the central support of $\omega$; then it follows that ${(e\otimes1_{L^\infty(\G)})\,\Delta(y)\in C_c(\G)\otimes C_c(\G)}$ and 
\eqnsp{(\mu\circ(\iota\otimes S))((e\otimes1_{L^\infty(\G)})\,\Delta(y))=\epsilon(y)\,e}
for all ${y\in C_c(\G)}$, which implies in conjunction with Lemma \ref{dqg-fin-supp-antipode-lemma-i} that 
\begin{align*}
(\chi_{\omega,2}\circ(\iota\otimes\Delta))(x\otimes y)&=\chi_{\omega,2}(x\otimes(e\otimes1_{L^\infty(\G)})\,\Delta(y))
\\&=(\omega\circ\mu\circ(\iota\otimes\mu)\circ(\iota\otimes\iota\otimes S))(x\otimes(e\otimes1_{L^\infty(\G)})\,\Delta(y))
\\&=(\omega\circ\mu)(x\otimes(\mu\circ(\iota\otimes S))((e\otimes1_{L^\infty(\G)})\,\Delta(y)))
\\&=(\omega\circ\mu)(x\otimes\epsilon(y)\,e)
\\&=\omega(xe)\,\epsilon(y)
\\&=(\omega\otimes\epsilon)(x\otimes y)
\end{align*}
for all ${x,y\in C_c(\G)}$. As ${C_c(\G)\otimes C_c(\G)}$ is ultraweakly dense in ${L^\infty(\G)\vnt L^\infty(\G)}$, it then follows by taking linear combinations and ultraweak limits that ${\chi_{\omega,2}\circ(\iota\otimes\Delta)=\omega\otimes\epsilon}$. 
\end{proof}

\begin{lemma}\label{dqg-fin-supp-antipode-lemma-iii}
Let $\G$ be amenable and ${\eta\in L^\infty(\G)^*}$ be a right invariant mean. If ${\omega\in\smash{L_f^1}(\G)}$, then 
\eqnsp{\eta\otimes\omega^\sharp=(\eta\otimes\chi_{\omega^*,1})\circ(\Delta\otimes\iota)}
\end{lemma}

\begin{proof}
Let ${e\in\Z(L^\infty(\G))}$ be the central support of $\omega$ and let ${f=S(e)}$; then 
\begin{align*}
(\eta\otimes\omega^\sharp)(x\otimes y)&=\eta(x)\,\omega^\sharp(y)
\\&=\omega^*(\eta(x)\,S(y))
\\&=\omega^*((\eta\otimes\iota)(\Delta(x))\,S(y))
\\&=\omega^*(e\,(\eta\otimes\iota)(\Delta(x))\,S(fy))
\\&=\chi_{\omega^*,1}((\eta\otimes\iota)(\Delta(x))\otimes y)
\\&=(\chi_{\omega^*,1}\circ(\eta\otimes\iota\otimes\iota)\circ(\Delta\otimes\iota))(x\otimes y)
\\&=((\eta\otimes\chi_{\omega^*,1})\circ(\Delta\otimes\iota))(x\otimes y)
\end{align*}
for all ${x\in L^\infty(\G)}$ and ${y\in C_c(\G)}$. Let ${X\in L^\infty(\G)\vnt L^\infty(\G)}$. Lemma \ref{central-supp-finite-dim} implies that ${e\in C_c(\G)}$, whence in turn ${f\in C_c(\G)}$; thus there exist ${x_1,\dots,x_n\in L^\infty(\G)}$ and ${y_1,\dots,y_n\in C_c(\G)}$ such that 
\eqnsp{(1_{L^\infty(\G)}\otimes f)\,X=\sum_{j\,=\,1}^n x_j\otimes y_j}
It therefore follows by Lemma \ref{omega-sharp-central-supp} that 
\begin{align*}
(\eta\otimes\omega^\sharp)(X)&=(\eta\otimes\omega^\sharp)((1_{L^\infty(\G)}\otimes f)\,X)
\\&=((\eta\otimes\chi_{\omega^*,1})\circ(\Delta\otimes\iota))((1_{L^\infty(\G)}\otimes f)\,X)
\\&=((\eta\otimes\chi_{\omega^*,1})\circ(\Delta\otimes\iota))(X)
\end{align*}
which proves the claimed equality. 
\end{proof}

We will also require two technical lemmas concerning tensor expansions. The first will allow us to expand products in such a way that medial tensor components are shifted onto distinct legs. 

\begin{lemma}\label{inverse-exp-tensor-exchange-lemma}
Let ${n\in\Nb}$. If ${X,Y\in\M\vnt L^\infty(\G)^{\vnt\msp(n+3)}}$ and ${\omega\in\smash{L_f^1}(\G)}$, then 
\eqnsp{(\iota^{\otimes\msp2}\otimes\chi_{\omega,1}\otimes\iota^{\otimes\msp n})(XY)=(\iota^{\otimes\msp2}\otimes\chi_{\omega,2}\otimes\iota^{\otimes\msp n})(X_{1235\,\cdots\,(n+5)}\hsp\msp Y_{1245\,\cdots\,(n+5)})}
\end{lemma}

\begin{proof}
Expanding tensor components demonstrates that 
\eqnsp{\chi_{\omega,1}((x\otimes y)(z\otimes w))=\chi_{\omega,2}((x\otimes y)_{13}\hsp(z\otimes w)_{23})}
for all ${x,y,z,w\in L^\infty(\G)}$. Tensoring on the left and right by ${\iota^{\otimes\msp2}}$ and ${\iota^{\otimes\msp n}}$ then yields the claimed equality in the case where $X$ and $Y$ are elementary tensors, and the general case follows by taking linear combinations and ultraweak limits. 
\end{proof}

The second will allow us to exchange iterated co-products with tensor expansions. 

\begin{lemma}\label{coproduct-tensor-exchange-lemma}
Let ${n\in\Nb}$. If ${X\in\M\vnt L^\infty(\G)}$, then 
\eqnsp{((\iota\otimes\Delta^{(n+1)})(X))_{124\,\cdots\,(n+4)}=(\iota^{\otimes\msp3}\otimes\Delta\otimes\iota^{\otimes\msp(n-1)})(((\iota\otimes\Delta^{(n)})(X))_{124\,\cdots\,(n+3)})}
\end{lemma}

\begin{proof}
Expanding tensor components demonstrates that 
\eqnsp{((\iota^{\otimes\msp(n+1)}\otimes\Delta)(x_0\otimes\cdots\otimes x_{n+1}))_{124\,\cdots\,(n+4)}=(\iota^{\otimes\msp(n+2)}\otimes\Delta)((x_0\otimes\cdots\otimes x_{n+1})_{124\,\cdots\,(n+3)})}
for all ${x_0\in\M}$ and ${x_1,\dots,x_{n+1}\in L^\infty(\G)}$; taking linear combinations and ultraweak limits, it then follows that 
\eqnsp{((\iota^{\otimes\msp(n+1)}\otimes\Delta)(Y))_{124\,\cdots\,(n+4)}=(\iota^{\otimes\msp(n+2)}\otimes\Delta)(Y_{124\,\cdots\,(n+3)})}
for all ${Y\in\M\vnt L^\infty(\G)^{\vnt\msp(n+1)}}$. If ${n=1}$, then taking ${Y=(\iota\otimes\Delta)(X)}$ yields the claimed equality. Suppose now that the claimed equality holds for some ${n\in\Nb}$; then 
\begin{align*}
((\iota\otimes\Delta^{(n+2)})(X))_{124\,\cdots\,(n+5)}&=((\iota^{\otimes\msp(n+2)}\otimes\Delta)((\iota\otimes\Delta^{(n+1)})(X)))_{124\,\cdots\,(n+5)}
\\&=(\iota^{\otimes\msp(n+3)}\otimes\Delta)(((\iota\otimes\Delta^{(n+1)})(X))_{124\,\cdots\,(n+4)})
\\&=(\iota^{\otimes\msp(n+3)}\otimes\Delta)((\iota^{\otimes\msp3}\otimes\Delta\otimes\iota^{\otimes\msp(n-1)})(((\iota\otimes\Delta^{(n)})(X))_{124\,\cdots\,(n+3)}))
\\&=(\iota^{\otimes\msp3}\otimes\Delta\otimes\iota^{\otimes\msp n})((\iota^{\otimes\msp(n+2)}\otimes\Delta)(((\iota\otimes\Delta^{(n)})(X))_{124\,\cdots\,(n+3)}))
\\&=(\iota^{\otimes\msp3}\otimes\Delta\otimes\iota^{\otimes\msp n})(((\iota^{\otimes\msp(n+1)}\otimes\Delta)((\iota\otimes\Delta^{(n)})(X)))_{124\,\cdots\,(n+4)})
\\&=(\iota^{\otimes\msp3}\otimes\Delta\otimes\iota^{\otimes\msp n})(((\iota\otimes\Delta^{(n+1)})(X))_{124\,\cdots\,(n+4)})
\end{align*}
The claim then follows by induction. 
\end{proof}

We are now prepared to prove our second factorization result. Throughout the remainder of this section, we assume that $\G$ is amenable and fix a right invariant mean ${\eta\in L^\infty(\G)^*}$. 

\begin{theorem}\label{dqg-answ-exp}
Let ${X,Y\in\M\vnt L^\infty(\G)}$. If ${Z\in\M\vnt L^\infty(\G)}$ is the element defined by 
\eqnsp{Z=(\iota\otimes\eta\otimes\iota)(Y_{12}^*\hsp\msp(\iota\otimes\Delta)(X))}
then 
\eqnsp{\Psi_Z(\theta^\sharp\star\omega)=(\iota\otimes\eta)((\iota\otimes\iota\otimes\theta)((\iota\otimes\Delta)(Y))^*\hsp\msp(\iota\otimes\iota\otimes\omega)((\iota\otimes\Delta)(X)))}
for all ${\omega,\theta\in\smash{L_\sharp^1}(\G)}$. 
\end{theorem}

\begin{proof}
It follows by applying (in order) Lemma \ref{dqg-fin-supp-antipode-lemma-iii}, Lemma \ref{inverse-exp-tensor-exchange-lemma} with ${n=1}$, Lemma \ref{coproduct-tensor-exchange-lemma} with ${n=2}$, Lemma \ref{dqg-fin-supp-antipode-lemma-ii}, and Lemma \ref{coproduct-tensor-exchange-lemma} with ${n=1}$ that 
\begin{align*}
\Psi_Z(\theta^\sharp\star\omega)&=((\iota\otimes\eta)\circ(\iota\otimes\iota\otimes(\theta^\sharp\star\omega)))(Y_{12}^*\hsp\msp(\iota\otimes\Delta)(X))
\\&=((\iota\otimes\eta)\circ(\iota\otimes\iota\otimes\theta^\sharp\otimes\omega))(Y_{12}^*\hsp\msp(\iota\otimes\Delta^{(2)})(X))
\\&=((\iota\otimes\omega)\circ(\iota\otimes\eta\otimes\theta^\sharp\otimes\iota))(Y_{12}^*\hsp\msp(\iota\otimes\Delta^{(2)})(X))
\\&=((\iota\otimes\omega)\circ(\iota\otimes((\eta\otimes\chi_{\theta^*,1})\circ(\Delta\otimes\iota))\otimes\iota))(Y_{12}^*\hsp\msp(\iota\otimes\Delta^{(2)})(X))
\\&=((\iota\otimes\omega)\circ(\iota\otimes\eta\otimes\chi_{\theta^*,1}\otimes\iota))(((\iota\otimes\Delta)(Y))_{123}^*\hsp\msp(\iota\otimes\Delta^{(3)})(X))
\\&=((\iota\otimes\eta\otimes\omega)\circ(\iota\otimes\iota\otimes\chi_{\theta^*,1}\otimes\iota))(((\iota\otimes\Delta)(Y))_{123}^*\hsp\msp(\iota\otimes\Delta^{(3)})(X))
\\&=((\iota\otimes\eta\otimes\omega)\circ(\iota\otimes\iota\otimes\chi_{\theta^*,2}\otimes\iota))(((\iota\otimes\Delta)(Y))_{123}^*\hsp\msp((\iota\otimes\Delta^{(3)})(X))_{12456})
\\&=((\iota\otimes\eta\otimes\omega)\circ(\iota\otimes\iota\otimes(\chi_{\theta^*,2}\circ(\iota\otimes\Delta))\otimes\iota))(((\iota\otimes\Delta)(Y))_{123}^*\hsp\msp((\iota\otimes\Delta^{(2)})(X))_{1245})
\\&=((\iota\otimes\eta\otimes\omega)\circ(\iota\otimes\iota\otimes\theta^*\otimes\epsilon\otimes\iota))(((\iota\otimes\Delta)(Y))_{123}^*\hsp\msp((\iota\otimes\Delta^{(2)})(X))_{1245})
\\&=((\iota\otimes\eta\otimes\omega)\circ(\iota\otimes\iota\otimes\theta^*\otimes((\epsilon\otimes\iota)\circ\Delta)))(((\iota\otimes\Delta)(Y))_{123}^*\hsp\msp((\iota\otimes\Delta)(X))_{124})
\\&=((\iota\otimes\eta\otimes\omega)\circ(\iota\otimes\iota\otimes\theta^*\otimes\iota))(((\iota\otimes\Delta)(Y))_{123}^*\hsp\msp((\iota\otimes\Delta)(X))_{124})
\\&=(\iota\otimes\eta)((\iota\otimes\iota\otimes\theta)((\iota\otimes\Delta)(Y))^*\hsp\msp(\iota\otimes\iota\otimes\omega)((\iota\otimes\Delta)(X)))
\end{align*}
for all ${\omega\in\smash{L_\sharp^1}(\G)}$ and ${\theta\in\smash{L_f^1}(\G)}$, which in conjunction with Lemma \ref{finite-supp-functional-dense-subspace} yields the claim. 
\end{proof}

As a consequence of the above theorem, we obtain the following analogue of \cite[Proposition 2.2]{dcot19} in the setting of amenable discrete quantum groups. 

\begin{cor}\label{dqg-twist-pos-def}
Let ${X\in\M\vnt L^\infty(\G)}$ be an isometry. The element ${Y\in\M\vnt L^\infty(\G)}$ defined by 
\eqnsp{Y=(\iota\otimes\eta\otimes\iota)(X_{12}^*\hsp\msp(\iota\otimes\Delta)(X))}
is a compression. 
\end{cor}

\begin{proof}
Fix ${n\in\Nb}$ and ${\omega\in M_n(\smash{L_\sharp^1}(\G))}$. Let ${Z\in M_n(\M\vnt L^\infty(\G))}$ be the element defined by 
\eqnsp{Z_{ij}=(\iota\otimes\iota\otimes\omega_{ij})((\iota\otimes\Delta)(X))}
It follows by Theorem \ref{dqg-answ-exp} that 
\eqnsp{\Psi_Y((\omega^\sharp)_{ik}\star\omega_{kj})=\Psi_Y((\omega_{ki})^\sharp\star\omega_{kj})=(\iota\otimes\eta)((Z_{ki})^*\hsp Z_{kj})=(\iota\otimes\eta)((Z^*)_{ik}\hsp Z_{kj})}
for ${i,j,k=1,\dots,n}$, whence 
\eqnsp{(\Psi_Y^{(n)}(\omega^\sharp\star\omega))_{ij}=\Psi_Y((\omega^\sharp\star\omega)_{ij})=(\iota\otimes\eta)((Z^*Z)_{ij})=((\iota\otimes\eta)^{(n)}(Z^*Z))_{ij}}
for ${i,j=1,\dots,n}$. As ${\iota\otimes\eta}$ is completely positive, this implies that 
\eqnsp{\Psi_Y^{(n)}(\omega^\sharp\star\omega)=(\iota\otimes\eta)^{(n)}(Z^*Z)\geq0}
and thus ${\Psi_Y}$ is completely positive. Furthermore, as ${(\iota\otimes\epsilon)\circ\Delta=\iota}$, it follows that 
\eqnsp{\Psi_Y(\epsilon)=(\iota\otimes\eta)((\iota\otimes\iota\otimes\epsilon)(X_{12}^*\hsp\msp(\iota\otimes\Delta)(X)))=(\iota\otimes\eta)(X^*X)=1_{\M\,\vnt\,L^\infty(\G)}}
and therefore ${\Psi_Y}$ is unital. The claim then follows from Proposition \ref{dqg-ucp-descends-compression}. 
\end{proof}

\section{Stability for approximate representations}\label{stability-section}

Throughout this section, we let $\G$ denote a locally compact quantum group and let $\M$ denote a von Neumann algebra. We will require the following notion of ``approximate unitarity''. 

\begin{defn}
Let $\A$ be a unital C*-algebra and ${\eps\geq0}$. An element ${u\in\A}$ is \emph{$\eps$-unitary} if 
\eqnsp{\norm{1_\A-u^*u}\leq\eps,\qquad\qquad\norm{1_\A-uu^*}\leq\eps}
\end{defn}

A key piece in the proof of Theorem \ref{cqg-adqg-stability} is the following lemma, which will allow us to upgrade a generic unitary approximate representation to an approximately unitary compression. 

\begin{lemma}\label{amenable-pos-def-approx}
Let ${\eps\geq0}$, ${\eta\in L^\infty(\G)^*}$ be a state, and ${X\in U(\M\vnt L^\infty(\G))}$ be an $\eps$-representation. The element ${Y\in\M\vnt L^\infty(\G)}$ defined by 
\eqnsp{Y=(\iota\otimes\eta\otimes\iota)(X_{12}^*\hsp\msp(\iota\otimes\Delta)(X))}
is ${\eps^2}$-unitary with ${\norm{Y}\leq1}$ and ${\norm{X-Y}\leq\eps}$. 
\end{lemma}

\begin{proof}
It is clear that ${\norm{Y}\leq1}$, and unitary invariance implies that 
\eqnsp{\norm{X_{12}^*\hsp\msp(\iota\otimes\Delta)(X)-X_{13}}=\norm{(\iota\otimes\Delta)(X)-X_{12}\hsp X_{13}}\leq\eps}
whence 
\begin{align*}
\norm{X-Y}&=\norm{(\iota\otimes\eta\otimes\iota)(X_{13}-X_{12}^*\hsp\msp(\iota\otimes\Delta)(X))}
\\&\leq\opnorm{\iota\otimes\eta\otimes\iota}\,\norm{X_{13}-X_{12}^*\hsp\msp(\iota\otimes\Delta)(X)}
\\&=\norm{X_{13}-X_{12}^*\hsp\msp(\iota\otimes\Delta)(X)}\leq\eps
\end{align*}
Furthermore, as 
\eqnsp{(\iota\otimes\eta\otimes\iota)(\abs{(\iota\otimes\Delta)(X)-X_{12}\hsp X_{13}}^2)=1_{\M\,\vnt\,L^\infty(\G)}-Y^*Y+\abs{X-Y}^2}
it follows that 
\eqnsp{0\leq1_{\M\,\vnt\,L^\infty(\G)}-Y^*Y\leq(\iota\otimes\eta\otimes\iota)(\abs{(\iota\otimes\Delta)(X)-X_{12}\hsp X_{13}}^2)}
wherefore 
\begin{align*}
\norm{1_{\M\,\vnt\,L^\infty(\G)}-Y^*Y}&\leq\norm{(\iota\otimes\eta\otimes\iota)(\abs{(\iota\otimes\Delta)(X)-X_{12}\hsp X_{13}}^2)}
\\&\leq\norm{\abs{(\iota\otimes\Delta)(X)-X_{12}\hsp X_{13}}^2}
\\&=\norm{(\iota\otimes\Delta)(X)-X_{12}\hsp X_{13}}^2\leq\eps^2
\end{align*}
Similarly, as 
\eqnsp{(\iota\otimes\eta\otimes\iota)(\abs{(X_{12}^*\hsp\msp(\iota\otimes\Delta)(X)-X_{13})^*}^2)=1_{\M\,\vnt\,L^\infty(\G)}-YY^*+\abs{(X-Y)^*}^2}
it follows that 
\begin{align*}
\norm{1_{\M\,\vnt\,L^\infty(\G)}-YY^*}&\leq\norm{(\iota\otimes\eta\otimes\iota)(\abs{(X_{12}^*\hsp\msp(\iota\otimes\Delta)(X)-X_{13})^*}^2)}
\\&\leq\norm{\abs{(X_{12}^*\hsp\msp(\iota\otimes\Delta)(X)-X_{13})^*}^2}
\\&=\norm{X_{12}^*\hsp\msp(\iota\otimes\Delta)(X)-X_{13}}^2\leq\eps^2
\end{align*}
and thus $Y$ is ${\eps^2}$-unitary. 
\end{proof}

The following lemma is an adaptation of \cite[Theorem 5.2]{dcot19} to our setting and constitutes the main technical result of this section. The key idea underlying the proof is taken from \cite[Theorem 3.1]{dcot19} (note however that our estimates differ significantly; cf. the proof of \cite[Theorem 2.6.5]{seguin22}). 

\begin{lemma}\label{amenable-pos-def-stability}
Let ${\eps\geq0}$ and $\G$ be amenable. If ${X\in\M\vnt L^\infty(\G)}$ is an $\eps$-unitary compression, then there exists a von Neumann algebra $\N$, two projections ${e,f\in\N}$, a *-isomorphism ${\sigma:\M\to e\N e}$, a partial isometry ${v\in f\N e}$, and a unitary representation ${U\in U(f\N f\vnt L^\infty(\G))}$ such that 
\eqnsp{\norm{(\sigma\otimes\iota)(X)-(v^*\otimes1_{L^\infty(\G)})\,U\hsp(v\otimes1_{L^\infty(\G)})}\leq4\eps}
and 
\eqnsp{\norm{e-v^*v}\leq\textstyle\frac{4}{3}\eps,\qquad\qquad\norm{f-vv^*}\leq4\eps}
\end{lemma}

\begin{proof}
As $X$ is a compression, it admits a Stinespring-type dilation ${(\N,e,\sigma,u,V)}$ such that $u$ is a partial isometry with ${u^*u=e}$. Let ${\eta\in L^\infty(\G)^*}$ be a right invariant mean and let ${q,y\in\N}$ be the elements defined by 
\eqnsp{q=uu^*,\qquad\qquad y=(\iota\otimes\eta)(V^*\,(q\otimes1_{L^\infty(\G)})\,V)}
Let ${f\in\N}$ be the projection defined by ${f=\chi_{[1/4,\msp1]}(y)}$. It follows by Lemma \ref{spectral-proj-ineq-ii} that 
\begin{align*}
\norm{e-u^*fu}&\leq\textstyle\frac{4}{3}\hsp\norm{e-u^*yu}
\\&=\textstyle\frac{4}{3}\hsp\norm{((\iota\otimes\eta)\circ(\sigma\otimes\iota))(1_\M\otimes1_{L^\infty(\G)}-X^*X)}
\\&\leq\textstyle\frac{4}{3}\hsp\norm{1_\M\otimes1_{L^\infty(\G)}-X^*X}\leq\frac{4}{3}\eps
\end{align*}
Let ${z=fu}$ and let ${z=v\msp\abs{z}}$ be the polar decomposition of $z$; then 
\eqnsp{\norm{e-\abs{z}^2}=\norm{e-z^*z}=\norm{e-u^*fu}\leq\textstyle\frac{4}{3}\eps}
As ${ze=z}$, it follows that ${z^*z\in e\N e}$, wherefore ${\abs{z}\in e\N e}$. As ${v^*v}$ is the support projection of ${\abs{z}}$ and ${\abs{z}\,(1_\N-e)=0}$, in turn 
\eqnsp{v\,(1_\N-e)=vv^*v\,(1_\N-e)=0}
and thus ${ve=v}$. Similar reasoning implies that ${\abs{z^*}\in f\N f}$ and ${v^*\,(1_\N-f)=0}$, and thus ${fv=v}$. It then follows by Lemma \ref{expand-ineq} that 
\eqnsp{\norm{e-v^*v}\leq\norm{e-v^*zz^*v}=\norm{e-\abs{z}\abs{z}^*}=\norm{e-\abs{z}^2}\leq\textstyle\frac{4}{3}\eps}
As 
\begin{align*}
y\otimes1_{L^\infty(\G)}&=(\iota\otimes\eta(\acdot)\,1_{L^\infty(\G)})(V^*\,(q\otimes1_{L^\infty(\G)})\,V)
\\&=(\iota\otimes((\eta\otimes\iota)\circ\Delta))(V^*\,(q\otimes1_{L^\infty(\G)})\,V)
\\&=((\iota\otimes\eta\otimes\iota)\circ(\iota\otimes\Delta))(V^*\,(q\otimes1_{L^\infty(\G)})\,V)
\\&=(\iota\otimes\eta\otimes\iota)(V_{13}^*\hsp V_{12}^*\,(q\otimes1_{L^\infty(\G)}\otimes1_{L^\infty(\G)})\,V_{12}\hsp V_{13})
\\&=V^*\,(y\otimes1_{L^\infty(\G)})\,V
\end{align*}
it follows that ${V\,(y\otimes1_{L^\infty(\G)})=(y\otimes1_{L^\infty(\G)})\,V}$, whence in turn ${V\,(f\otimes1_{L^\infty(\G)})=(f\otimes1_{L^\infty(\G)})\,V}$. Let ${U\in U(f\N f\vnt L^\infty(\G))}$ be the unitary representation defined by 
\eqnsp{U=(f\otimes1_{L^\infty(\G)})\,V\hsp(f\otimes1_{L^\infty(\G)})}
As 
\begin{align*}
\norm{(\sigma\otimes\iota)(X)-(z^*\otimes1_{L^\infty(\G)})\,U\hsp(z\otimes1_{L^\infty(\G)})}&=\norm{(\sigma\otimes\iota)(X)-(u^*\otimes1_{L^\infty(\G)})\,U\hsp(u\otimes1_{L^\infty(\G)})}
\\&=\norm{(u^*\otimes1_{L^\infty(\G)})(V-U)(u\otimes1_{L^\infty(\G)})}
\\&=\norm{(u^*f^\perp\otimes1_{L^\infty(\G)})\,V\hsp(f^\perp u\otimes1_{L^\infty(\G)})}
\\&\leq\norm{u^*f^\perp\otimes1_{L^\infty(\G)}}\,\norm{V}\,\norm{f^\perp u\otimes1_{L^\infty(\G)}}
\\&=\norm{u^*f^\perp\otimes1_{L^\infty(\G)}}\,\norm{f^\perp u\otimes1_{L^\infty(\G)}}
\\&=\norm{u^*f^\perp}\,\norm{f^\perp u}
\\&=\norm{u^*f^\perp u}\leq\textstyle\frac{4}{3}\eps
\end{align*}
it follows by Lemma \ref{square-ineq} that 
\begin{align*}
\norm{(\sigma\otimes\iota)(X)-(v^*\otimes1_{L^\infty(\G)})\,U\hsp(v\otimes1_{L^\infty(\G)})}&\leq\norm{(z^*\otimes1_{L^\infty(\G)})\,U\hsp((z-v)\otimes1_{L^\infty(\G)})}
\\&\quad+\norm{((z-v)^*\otimes1_{L^\infty(\G)})\,U\hsp(v\otimes1_{L^\infty(\G)})}+\textstyle\frac{4}{3}\eps
\\&=\norm{(z^*\otimes1_{L^\infty(\G)})\,U\hsp(v\,(\abs{z}-e)\otimes1_{L^\infty(\G)})}
\\&\quad+\norm{((v\,(\abs{z}-e))^*\otimes1_{L^\infty(\G)})\,U\hsp(v\otimes1_{L^\infty(\G)})}+\textstyle\frac{4}{3}\eps
\\&\leq\norm{(z^*\otimes1_{L^\infty(\G)})\,U\hsp(v\otimes1_{L^\infty(\G)})}\,\norm{\abs{z}-e}
\\&\quad+\norm{\abs{z}-e}\,\norm{(v^*\otimes1_{L^\infty(\G)})\,U\hsp(v\otimes1_{L^\infty(\G)})}+\textstyle\frac{4}{3}\eps
\\&\leq2\,\norm{e-\abs{z}}+\textstyle\frac{4}{3}\eps
\\&\leq2\,\norm{e-\abs{z}^2}+\textstyle\frac{4}{3}\eps\leq4\eps
\end{align*}
Furthermore, Lemma \ref{expand-ineq} and Lemma \ref{spectral-proj-ineq} imply that 
\begin{align*}
\norm{f-vv^*}&\leq\norm{f-v\msp\abs{z}\abs{z}\msp v^*}
\\&=\norm{f-zz^*}
\\&=\norm{f\,(1_\N-q)\,f}
\\&=\norm{(1_\N-q)\,f\,(1_\N-q)}
\\&\leq4\,\norm{(1_\N-q)\,y\,(1_\N-q)}
\\&=4\,\norm{(\iota\otimes\eta)((q^\perp\otimes1_{L^\infty(\G)})\,V^*\hsp(q\otimes1_{L^\infty(\G)})\,V\hsp(q^\perp\otimes1_{L^\infty(\G)}))}
\\&\leq4\,\norm{(q^\perp\otimes1_{L^\infty(\G)})\,V^*\hsp(q\otimes1_{L^\infty(\G)})\,V\hsp(q^\perp\otimes1_{L^\infty(\G)})}
\\&=4\,\norm{(u^*\otimes1_{L^\infty(\G)})\,V\hsp(q^\perp\otimes1_{L^\infty(\G)})\,V^*\hsp(u\otimes1_{L^\infty(\G)})}
\\&=4\,\norm{e\otimes1_{L^\infty(\G)}-(u^*\otimes1_{L^\infty(\G)})\,V\hsp(q\otimes1_{L^\infty(\G)})\,V^*\hsp(u\otimes1_{L^\infty(\G)})}
\\&=4\,\norm{(\sigma\otimes\iota)(1_\M\otimes1_{L^\infty(\G)}-XX^*)}
\\&=4\,\norm{1_\M\otimes1_{L^\infty(\G)}-XX^*}\leq4\eps
\end{align*}
which proves the claim. 
\end{proof}

We are now ready to prove the stability results from the introduction. We begin with the proof of Theorem \ref{amenable-pos-def-stability-thm}, which is an easy consequence of the preceding lemma. 

\begin{proof}[Proof of Theorem \ref{amenable-pos-def-stability-thm}]
By Lemma \ref{amenable-pos-def-stability}, there exists a von Neumann algebra $\N$, projections ${e,f\in\N}$, a *-isomorphism ${\sigma:\M\to e\N e}$, a partial isometry ${v\in\N}$ such that ${v^*v\leq e}$ and ${vv^*\leq f}$, and a unitary representation ${U\in U(f\N f\vnt L^\infty(\G))}$ such that 
\eqnsp{\norm{(\sigma\otimes\iota)(X)-(v^*\otimes1_{L^\infty(\G)})\,U\hsp(v\otimes1_{L^\infty(\G)})}\leq4\eps}
and 
\eqnsp{\norm{e-v^*v}<1,\qquad\qquad\norm{f-vv^*}<1}
The latter two inequalities imply that ${v^*v=e}$ and ${vv^*=f}$, and thus ${(v^*\otimes1_{L^\infty(\G)})\,U\hsp(v\otimes1_{L^\infty(\G)})}$ is a unitary representation in ${e\N e\vnt L^\infty(\G)}$. 
\end{proof}

With Theorem \ref{amenable-pos-def-stability-thm} in hand, we at last have all the necessary pieces to prove Theorem \ref{cqg-adqg-stability}. 

\begin{proof}[Proof of Theorem \ref{cqg-adqg-stability}]
Let ${\eta\in L^\infty(\G)^*}$ be a right invariant mean and let ${Y\in\M\vnt L^\infty(\G)}$ be the element defined by 
\eqnsp{Y=(\iota\otimes\eta\otimes\iota)(X_{12}^*\hsp\msp(\iota\otimes\Delta)(X))}
If $\G$ is compact, then ${\eta=\psi}$ by the uniqueness of the right Haar weight up to rescaling, and thus it follows by Corollary \ref{cqg-twist-pos-def} that $Y$ is a compression; if $\G$ is instead an amenable discrete quantum group, then it follows by Corollary \ref{dqg-twist-pos-def} that $Y$ is a compression. In either case, the claim then follows from Lemma \ref{amenable-pos-def-approx} in conjunction with Theorem \ref{amenable-pos-def-stability-thm}. 
\end{proof}

\end{document}